\documentclass[11pt,twoside,letterpaper]{article} 
\usepackage{times,fancyhdr}   
\usepackage{amsfonts}                     
\usepackage{amsthm}                      
\usepackage{amsmath}                   
\usepackage{amssymb}  

\numberwithin{equation}{section}


\makeatletter 
\def\cleardoublepage{\clearpage\if@twoside \ifodd\c@page\else%
    \hbox{}%
    \thispagestyle{empty}%
    \newpage%
    \if@twocolumn\hbox{}\newpage\fi\fi\fi} 
\makeatother 

\def \C{\mathbb{C}}

\def \Q{\mathbb{Q}}

\def \Z{\mathbb{Z}}

\newtheorem{thm}{Theorem}[section]

\newtheorem{lemma}{Lemma}[section]

\begin{document}
\title{
{\begin{flushleft}
\vskip 0.45in
{\normalsize\bfseries\textit{ }}
\end{flushleft}
\vskip 0.45in
\bfseries\scshape Algebraic independence results for values of theta-constants, II}}

\thispagestyle{fancy}
\fancyhead{}
\fancyhead[L]{In: Book Title \\ 
Editor: Editor Name, pp. {\thepage-\pageref{lastpage-01}}} 
\fancyhead[R]{ISBN 0000000000  \\
\copyright~2007 Nova Science Publishers, Inc.}
\fancyfoot{}
\renewcommand{\headrulewidth}{0pt}

\author{\bfseries\itshape Carsten Elsner\thanks{Fachhochschule f{\"u}r die Wirtschaft, University of Applied Sciences, Freundallee 15,
D-30173 Hannover, Germany \newline e-mail: carsten.elsner@fhdw.de},
Yohei Tachiya\thanks{Hirosaki University, Graduate School of Science and Technology, 
Hirosaki 036-8561, Japan \newline e-mail: tachiya@hirosaki-u.ac.jp}}

\date{}
\maketitle
\thispagestyle{empty}
\setcounter{page}{1}

\begin{center}
Dedicated to Professor Iekata Shiokawa on the occasion of his 75th birthday
\end{center}
\[\]
\[\]
\begin{abstract}
Let $\theta_3(\tau)=1+2\sum_{\nu=1}^{\infty} q^{\nu^2}$ with $q=e^{i\pi \tau}$ denote the Thetanullwert of the Jacobi 
theta function 
\[\theta(z|\tau) \,=\,\sum_{\nu=-\infty}^{\infty} e^{\pi i\nu^2\tau + 2\pi i\nu z} \,.\]
Moreover, let $\theta_2(\tau)=2\sum_{\nu=0}^{\infty} q^{{(\nu+1/2)}^2}$ and $\theta_4(\tau)=1+2\sum_{\nu=1}^{\infty} {(-1)}^{\nu}q^{\nu^2}$.
For algebraic numbers $q$ with $0<|q|<1$ and for any $j\in \{ 2,3,4\}$ we prove the algebraic independence over ${\Q}$ of the numbers 
$\theta_j(n\tau)$ and $\theta_j(\tau)$ for all odd integers $n\geq 3$. Assuming the same conditions on $q$ and $\tau$ as above,
we obtain sufficient conditions by use of a criterion involving resultants in order to decide on the algebraic independence over ${\Q}$ of 
$\theta_j(2m\tau)$ and $\theta_j(\tau)$ $(j=2,3,4)$ and of $\theta_3(4m\tau)$ and $\theta_3(\tau)$ with odd positive integers $m$. 
In particular, we prove the algebraic independence of $\theta_3(n\tau)$ and $\theta_3(\tau)$ for even integers $n$ with $2\leq n\leq 22$. 
The paper continues the work of the first-mentioned author, who already proved the algebraic independence of $\theta_3(2^m\tau)$ and 
$\theta_3(\tau)$ for $m=1,2,\dots$. 
\end{abstract}
 
\vspace{.08in} \noindent \textbf{Keywords:} Algebraic independence, Theta-constants, Nesterenko's theorem, 
Independence criterion,  \newline \hspace*{52pt} Modular equations 


\pagestyle{fancy}  
\fancyhead{}
\fancyhead[EC]{Carsten Elsner}
\fancyhead[EL,OR]{\thepage}
\fancyhead[OC]{Algebraic independence results for values of theta-constants, II}
\fancyfoot{}
\renewcommand\headrulewidth{0.5pt} 

\newpage
\section{Introduction and statement of results} \label{Sec1}   
Let \(\tau \) be a complex variable in the complex upper 
half-plane \(\Im (\tau) >0 \). 
The series
\[\theta_2(\tau) = 2\sum_{\nu=0}^{\infty} q^{{(\nu+1/2)}^2} \,,\qquad \theta_3(\tau) = 1+2\sum_{\nu=1}^{\infty} 
q^{\nu^2} \,,\qquad \theta_4(\tau) = 1+2\sum_{\nu=1}^{\infty} {(-1)}^\nu q^{\nu^2} \]
are known as theta-constants or Thetanullwerte, where \(q=e^{\pi i\tau} \). In particular,  $\theta_3(\tau)$ 
is the Thetanullwert of the Jacobi theta function $\theta(z|\tau) \,=\,\sum_{\nu=-\infty}^{\infty} e^{\pi i\nu^2\tau + 2\pi i\nu z}$. 
For an extensive discussion of theta-functions and theta-constants we refer the reader to \cite{Hurwitz}, \cite{Lang}, and \cite{Lawden}.
Recently, the first-named author has proven the following result.\\

{\bf Theorem~A.} \cite[Theorem~1.1]{Elsner3}
{\em Let $q$ be an algebraic number with \(q=e^{\pi i \tau} \) and \(\Im (\tau) >0 \). Let  $m\geq1$ be an integer. Then, the two numbers $\theta_3(2^m\tau)$ and
$\theta_3(\tau)$ are algebraically independent over ${\Q}$ as well as the two numbers $\theta_3(n\tau)$ and $\theta_3(\tau)$ for $n=3,5,6,7,9,10,11,12$\/}.\\

The first basic tool in proving such algebraic independence results are integer polynomials in two variables $X,Y$, which vanish at certain points $X=X_0$ and $Y=Y_0$
given by values of rational functions of theta-constants. For instance, for $n=2^m$ 
$(m\geq3)$ we consider the polynomial
\[P_n(X,Y) \,=\, {\big( nX - {(1+Y)}^2 \big)}^{2^{m-2}} + YU_n\big( X,{(1+Y)}^2,Y \big) \]
where $U_n(t_1,t_2,t_3)\in {\Q}[t_1,t_2,t_3]$ is a polynomial satisfying 
\[U_n\Big( \,\frac{1}{n},1,0\,\Big) \,=\, -2^{2^{m-1}-1} \,.\]
Moreover, we have
\[P_n\Big( \,\frac{\theta_3^2(n\tau)}{\theta_3^2(\tau)},\,\frac{\theta_4(\tau)}{\theta_3(\tau)} \,\Big) \,=\, 0\,.\]
The second tool is an algebraic independence criterion (see Lemma~\ref{Lem1} below), from which the algebraic independence of $\theta_3(n\tau)$ and $\theta_3(\tau)$
over ${\Q}$ can be obtained by proving that the resultant
$$
\mbox{Res}_X\,\Big( \,P_n(X,Y),\,\frac{\partial}{\partial Y}P_n(X,Y)\,\Big) \,\in \, {\Z}[Y] 
$$
does not vanish identically (see \cite[Theorem~4.1]{Elsner3}). This can be seen as follows. We have the identities
\begin{eqnarray*}
P_n(X,0) &=& {\big( 2^mX-1 \big)}^{2^{m-2}} \,,\\
\frac{\partial P_n}{\partial Y} (X,0) &=& -2^{m-1}{\big( 2^mX-1 \big)}^{2^{m-2}-1} + U_n(X,1,0) \,,
\end{eqnarray*}
from which on the one hand we deduce that $P_n(X,0)$ has a $2^{m-2}$-fold root at $X_1=1/2^m$. On the other hand one has
\[\frac{\partial P_n}{\partial Y} (X_1,0) \,=\, U_n\Big( \,\frac{1}{n},1,0\,\Big) \,=\, -2^{2^{m-1}-1} \,\not= \, 0\,.\]
Hence, for $Y=0$ the polynomials $P_n(X,Y)$ and $\partial P_n(X,Y)/\partial Y$ 
have no common root. Therefore, the above resultant with respect to $X$ does not vanish identically, which gives the desired result. 
We state the polynomials $P_{2^m}(X,Y)$ for $m=1,2,3,4$ explicitly. The algorithm to compute these polynomials recursivly is given by Lemma~3.1 in \cite{Elsner3}. 
\begin{eqnarray*}
P_2 &=& 2X-Y^2-1 \,,\\ 
P_4 &=& 4X-{(1+Y)}^2 \,,\\
P_8 &=& 64X^2 - 16{(1+Y)}^2X + {(1-Y)}^4 \,,\\
P_{16} &=& 65536X^4 - 16384{(1+Y)}^2X^3 + 512(3Y^4+4Y^3+18Y^2+4Y+3)X^2 \\
&& -\,64{(1+Y)}^2(Y^4+28Y^3+6Y^2+28Y+1)X + {(1-Y)}^8 \,,
\end{eqnarray*}

Let \(n\geq 3 \) denote an odd positive integer. Set
\[h_j(\tau) := n^2\frac{\theta_j^4(n\tau)}{\theta_j^4(\tau)} \quad (j=2,3,4)\,,\quad 
\lambda = \lambda(\tau) := \frac{\theta_2^4(\tau)}{\theta_3^4(\tau)} \,,
\quad \psi (n) := n\prod_{p|n} \Big( 1+\frac{1}{p} \Big) \,,\]
where \(p\) runs through all primes dividing \(n\). Yu.V. Nesterenko \cite{Nes3} proved the existence of integer polynomials $P_n(X,Y)\in
{\Z}[X,Y]$ such that $P_n\big(h_j(\tau),R_j(\lambda(\tau))\big)=0$ holds for $j=2,3,4$, odd integers $n\geq 3$,
and a suitable rational function $R_2,R_3$, or $R_4$, respectively.\\

{\bf Theorem~B.} \cite[Theorem~1, Corollary~3]{Nes3} 
{\em For any odd integer \(n\geq 3 \) there exists a polynomial \(P_n(X,Y) \in {\Z}[X,Y] \), \(\deg_X P_n = \psi (n) \),
such that\/}
\begin{eqnarray}
P_n\Big( \,h_2(\tau),16\frac{\lambda(\tau)-1}{\lambda(\tau)} \,\Big) &\,=\,& 0\,, \label{eq:101}\\
P_n\big( h_3(\tau),16\lambda(\tau) \big) &\,=\,& 0 \,,
\label{eq:102}\\
P_n\Big( \,h_4(\tau),16\frac{\lambda(\tau)}{\lambda(\tau)-1} \,\Big) &\,=\,& 0 \,.\label{eq:103}
\end{eqnarray}
The polynomials $P_3,P_5,P_7,P_9$, and $P_{11}$ are listed in the appendix of \cite{Elsner3}. $P_3$ and $P_5$ are already given in \cite{Nes3},
$P_7,P_9$, and $P_{11}$ are the results of computer-assisted computations of the first-named author. \\

In this paper we focus on the problem to decide on the algebraic independence of $\theta_j(n\tau)$ and
$\theta_j(\tau)$ $(j=2,3,4)$ over ${\Q}$ for algebraic numbers $q$, 
where $n\geq 3$ is an odd integer or $n=2m,4m$ with odd positive integers $m$.
The above Theorem~B will be used in Section \ref{Sec2}.\\

In the following theorems, 
the number \(q=e^{\pi i \tau} \) is an algebraic number with \(\Im (\tau) >0 \).
\begin{thm}\label{Thm1}
Let $n\geq3$ be an odd integer. Then, the numbers in each of the sets 
$$
\{\theta_2(n\tau),\theta_2(\tau)\},\quad\{\theta_3(n\tau),\theta_3(\tau)\},\quad 
\{\theta_4(n\tau),\theta_4(\tau)\}
$$
are algebraically independent over $\Bbb{Q}$.
\end{thm}
In order to prove this theorem we first shall show that for an algebraic number $q$ with $0<|q|<1$ the numbers $h_2(\tau)$, $h_3(\tau)$, and $h_4(\tau)$ are transcendental (Lemma~\ref{Lem7}).
This interim result already shows that the two numbers 
$\theta_j(n\tau)$ and $\theta_j(\tau)$ $(j=2,3,4)$ are homogeneously algebraically independent over ${\Q}$. 

On the other hand, it has not been shown that Theorem \ref{Thm1} holds 
for arbitrary even integers $n$. However we can prove it 
for small even integers $n$ by checking the non-vanishing of a Jacobian determinant (Lemma~\ref{Lem1}), which is hard to decide when the 
involved polynomials are not given explicitly. 
\begin{thm}
For $n=2,4,6$, 
the numbers $\theta_2(n\tau)$ and $\theta_2(\tau)$ are algebraically independent over ${\Q}$.
\label{Thm2}
\end{thm}
\begin{thm}
For $n=2,4,6,8,10$, 
the numbers $\theta_4(n\tau)$ and $\theta_4(\tau)$ are algebraically independent over ${\Q}$.
\label{Thm4}
\end{thm}
\begin{thm}
Let  $2\leq n\leq 22$ be an even integer. Then,
the numbers $\theta_3(n\tau)$ and $\theta_3(\tau)$ are algebraically independent over ${\Q}$.
\label{Thm3}
\end{thm}
\section{Auxiliary results} \label{Sec2}
In this section, we prepare some lemmas to prove theorems.
\begin{lemma} \cite[Lemma~4]{Elsner2}
Let $q$ be an algebraic number with \(q=e^{\pi i \tau} \) and \(\Im (\tau) >0 \). 
Then, any two numbers in the set  
in each of the sets  
\[\big\{ \theta_2(\tau),\theta_3(\tau),\theta_4(\tau) \big\}  \]
are algebraically independent over \({\Q} \).
\label{Lem2}
\end{lemma}
This result can be derived from Yu.V. Nesterenko's theorem \cite{Nes1} on the algebraic independence of the values
$P(q),Q(q),R(q)$ of the Ramanujan functions $P,Q,R$ at a nonvanishing algebraic point $q$. It should be noticed that the three numbers $\theta_2(\tau)$, $\theta_3(\tau)$, and $\theta_4(\tau)$ are algebraically dependent over $\Bbb{Q}$, since the identity 
\begin{equation}\label{eq:405}
\theta_3^4(\tau)=\theta_2^4(\tau)+\theta_4^4(\tau)
\end{equation}
holds for any $\tau\in\Bbb{C}$ with \(\Im (\tau) >0 \).\\

In what follows, we distinguish two cases based on the parity of $n$. 
\subsection{The case where $n$ is odd}
The following subsequent Lemmas \ref{Lem6}, \ref{Lem7}, and \ref{Lem71} are needed to prove Theorem \ref{Thm1}. 
Let $n\geq3$ be a fixed odd integer and $\tau \in {\C}$ with $\Im (\tau)>0$. 
From
Theorem~B we know that there exists a nonzero polynomial $P_n(X,Y)\in {\Z}[X,Y]$ 
with $\deg_XP_n=\psi(n)$ such that $P_n(X_0,Y_0)$ vanishes for
$$
X_0 \,:=\, h_3(\tau) \,=\, n^2\frac{\theta_3^4(n\tau)}{\theta_3^4(\tau)}\,,\qquad Y_0 \,:=\, 16\lambda(\tau) \,=\, 16\frac{\theta_2^4(\tau)}{\theta_3^4(\tau)} \,.
$$
Let $N:=\deg_Y P_n(X,Y)$. The polynomials $Q_j(X)\in {\Z}[X]$ $(j=0,1,\dots ,N)$ are given by
\begin{equation}
P_n(X,Y) \,=\, \sum_{j=0}^N Q_j(X)Y^j \,.
\label{2.10}
\end{equation}
\begin{lemma}
For any complex number $\alpha$, there exists a subscript $j=j(\alpha)$ such that $Q_j(\alpha)\not= 0$. 
\label{Lem6}
\end{lemma}
{\em Proof of Lemma~\ref{Lem6}.\/} 
\,Suppose on the contrary that there exists an $\alpha\in\Bbb{C}$ such that $Q_j(\alpha)=0$ for all $j=0,1,\dots ,N$.
It follows from (\ref{2.10}) that there exists a polynomial $R_n(X,Y)\in {\C}[X,Y]$ satisfying
\begin{equation}
P_n(X,Y) \,=\, (X-\alpha) R_n(X,Y) \,.
\label{H120}
\end{equation}
In accordance with formula (5) in \cite{Nes3} we define for any $\tau \in {\C}$ with $\Im (\tau)>0$ the numbers
\[x_{\nu}(\tau) \,:=\, u^2\frac{\theta_3^4\big( \frac{u\tau +2v}{w}\big)}{\theta_3^4(\tau)} \qquad \big( \nu =1,2,\dots,\psi(n)\big) \,,\]
where the nonnegative integers $u,v,w$ are given by \cite[Lemma~1]{Nes3}. These integers depend on $n$ and $\nu$ and satisfy the three conditions
\begin{equation}
(u,v,w) \,=\, 1 \,,\quad uw \,=\, n \,,\quad 0\leq v<w \,.
\label{H130}
\end{equation}
Substituting $Y=16\lambda (\tau)$ into (\ref{H120}), we have by \cite[Corollary~1]{Nes3}
\begin{equation}
\prod_{\nu =1}^{\psi(n)} \big( X - x_{\nu}(\tau)\big) \,=\, P_n\big( X,16\lambda(\tau)\big) \,=\, (X-\alpha)R_n\big( X,16\lambda(\tau) \big) \,.
\label{H140}
\end{equation}
Next, by substituting $X=\alpha$ into (\ref{H140}), we obtain
\begin{equation}
\prod_{\nu =1}^{\psi(n)} \big( \alpha - x_{\nu}(\tau)\big) \,=\, 0 
\label{H150}
\end{equation}
for any $\tau \in {\C}$ with $\Im (\tau)>0$. Let $a_k:=nki$ $(k=1,2,\dots)$ be a sequence of complex numbers on the imaginary axis. Then we get by (\ref{H150})
\[\prod_{\nu =1}^{\psi(n)} \big( \alpha - x_{\nu}(a_k)\big) \,=\, 0 \qquad (k=1,2,\dots) \,.\]
Hence, by the pigeonhole principle, there is a subscript $\nu_0$ with $1\leq \nu_0 \leq \psi(n)$ such that
\begin{equation}
u^2\frac{\theta_3^4\big( \frac{ub_k +2v}{w}\big)}{\theta_3^4(b_k)} \,=\, x_{\nu_0}(b_k) \,=\, \alpha 
\label{H160}
\end{equation}
holds for some subsequence ${\{ b_k\}}_{k\geq 1}$ of ${\{ a_k\}}_{k\geq 1}$. The integers $u,v,w$ in (\ref{H160}) depend on $n,\nu_0$ and satisfy the condions in
(\ref{H130}). Since $\sqrt[4]{\alpha/u^2}$ takes four complex values, we see in the same way by applying the pigeonhole principle that there exists a complex
number $\beta$ with $\beta^4 =\alpha/u^2$ such that
\begin{equation}
\frac{\theta_3\big( \frac{uc_k +2v}{w}\big)}{\theta_3(c_k)} \,=\, \beta 
\label{H170}
\end{equation}
holds for some subsequence ${\{ c_k\}}_{k\geq 1}$ of ${\{ b_k\}}_{k\geq 1}$. 
Let $c_k := nt_k i$, where $t_k$ $(k\geq1)$ are positive integers with $t_1<t_2<\dots$.
Then, by substituting $\tau =c_k$ into
$q=e^{\pi i\tau}=e^{-\pi t_kn}$, we obtain
\begin{equation}
\theta_3(c_k) \,=\, 1+2\sum_{m=1}^{\infty} {\big( e^{-\pi t_k} \big)}^{nm^2} \,.
\label{H180}
\end{equation}
For $\tau = \frac{uc_k +2v}{w}$ it follows with $\xi_w :=e^{\frac{2\pi i}{w}}$ that
\[q \,=\, e^{\pi i\frac{unt_ki+2v}{w}} \,=\, {\big( e^{\frac{2\pi i}{w}}\big)}^v \cdot e^{-\pi u^2t_k} \,=\, \xi_w^v \cdot e^{-\pi u^2t_k} \,,\]
which yields
\begin{equation}
\theta_3\Big( \,\frac{uc_k+2v}{w}\,\Big) \,=\, 1+2\sum_{m=1}^{\infty} {\xi}_w^{vm^2} {\big( e^{-\pi t_k} \big)}^{u^2m^2} \,. 
\label{H190}
\end{equation}
Next, two complex functions $f(z)$ and $g(z)$ are defined by their power series, namely
\begin{eqnarray*}
f(z) &:=& 1+2\sum_{m=1}^{\infty} z^{nm^2} \,,\\
g(z) &:=& 1+2\sum_{m=1}^{\infty} {\xi}_w^{vm^2}z^{u^2m^2} \,.
\end{eqnarray*}
Setting $\eta_k:=e^{-\pi t_k}$ $(k=1,2,\dots)$, it follows from (\ref{H170}) to (\ref{H190}) that
\[\beta f(\eta_k) \,=\, g(\eta_k) \qquad (k=1,2,\dots) \,.\]
The sequence ${\{ \eta_k\}}_{k\geq 1}$ tends to zero, such that we may apply the identity theorem for power series. We obtain
\[\beta f(z) \,=\, g(z) \qquad (|z|<1) \,.\]
Comparing the first and second nonvanishing coefficients of the series, it follows that $\beta =1$, and, by applying $uw=n$ in (\ref{H130}), $u^2=n$, $u=w=\sqrt{n}>1$,
$\xi_w^v=1$. Since $\xi_w \not=1$, we conclude from $\xi_w^v=1$ and $0\leq v<w$ in (\ref{H130}) that $v=0$. Finally, we deduce that $(u,v,w)=(u,0,u)=u>1$, which 
contradicts the arithmetic condition $(u,v,w)=1$ in (\ref{H130}). This completes the proof of Lemma~\ref{Lem6}. \hfill $\Box$
\par\medskip
\begin{lemma}
Let $q$ be an algebraic number with \(q=e^{\pi i \tau} \) and \(\Im (\tau) >0 \). 
Then the three numbers $h_2(\tau)$, $h_3(\tau)$, and $h_4(\tau)$ are transcendental.
\label{Lem7}
\end{lemma}
{\em Proof of Lemma~\ref{Lem7}.\/}  
\,We suppose on the contrary that $h_3(\tau)$ is an algebraic number. Then it follows from (\ref{2.10}) and from Lemma~\ref{Lem6}
that
\[F(Y) \,:=\, P_n\big( h_3(\tau),Y \big) \,=\, \sum_{j=0}^N Q_j\big( h_3(\tau) \big) Y^j \]
is a nonzero polynomial with algebraic coefficients. 
Hence, the identity (\ref{eq:102}) yields $F\big( 16\lambda(\tau) \big)=0$, thus showing that 
$\lambda(\tau)=\frac{\theta_2^4(\tau)}{\theta_3^4(\tau)}$ is an algebraic number. But, by Lemma~\ref{Lem2}, the numbers $\theta_2(\tau)$ and $\theta_3(\tau)$ are
algebraically independent over ${\Q}$, a contradiction. Thus, $h_3(\tau)$ is transcendental. 
Similarly, the transcendence of $h_2(\tau)$ and $h_4(\tau)$ follows from 
the identities (\ref{eq:101}) and (\ref{eq:103}). \hfill $\Box$  
\begin{lemma}\label{Lem71}
Let $\alpha_1,\alpha_2\in\Bbb{C}$ be algebraically independent over $\Bbb{Q}$ 
and let $\beta_1,\beta_2\in\Bbb{C}$ $(\beta_2\neq0)$ such that $\beta_1/\beta_2$ is transcendental. 
Suppose that there exist nonzero polynomials  
$P(X,Y),Q(X,Y)\in\Bbb{Q}[X,Y]$ such that 
\begin{equation}\label{eq:10}
P(\alpha_1/\alpha_2,\beta_1/\beta_2)=0.
\end{equation}
and 
\begin{equation}\label{eq:20}
Q(\beta_1,\beta_2)=\alpha_2.
\end{equation}
Then $\beta_1$ and $\beta_2$ are algebraically independent over $\Bbb{Q}$.
\end{lemma}
{\em Proof of Lemma~\ref{Lem71}.\/}
\,Define the fields $F:=\mathbb{Q}(\beta_1,\beta_2)$ and $E:=F(\alpha_1,\alpha_2)$. 
We first prove that the field extension $E/F$ is algebraic. 
Since $\alpha_2\in F$ by (\ref{eq:20}), we only have to show that $\alpha_1$ is algebraic over $F$. 
Let 
$$
P(X,Y):=\sum_{j=0}^{\ell}R_j(Y)X^j,\qquad R_j(Y)\in\mathbb{Q}[Y],\qquad R_\ell(Y)\not\equiv0,
$$
and define 
$$
f(X):=\sum_{j=0}^{\ell}R_j(\beta_1/\beta_2)(X/\alpha_2)^j\in F[X],
$$
where $f(X)$ is a nonzero polynomial, since 
$R_\ell(\beta_1/\beta_2)\neq0$ follows from the transcendence of $\beta_1/\beta_2$. 
By (\ref{eq:10}), we have $f(\alpha_1)=P(\alpha_1/\alpha_2,\beta_1/\beta_2)=0$, which implies that $\alpha_1$ is algebraic over $F$. 

Thus we get  
$$
\begin{array}{ll}
{\rm trans.}\deg F/\mathbb{Q}={\rm trans.}\deg E/F+
{\rm trans.}\deg F/{\mathbb{Q}}=
{\rm trans.}\deg E/{\mathbb{Q}}\geq2,
\end{array}
$$
where we used the algebraic independence hypothesis on $\alpha_1$ and $\alpha_2$. 
On the other hand, ${\rm trans.}\deg F/\mathbb{Q}\leq2$ is trivial. 
Therefore we obtain 
$${\rm trans.}\deg F/\mathbb{Q}=2,$$
which gives the desired result. 
\hfill \qed
\subsection{The case where $n$ is even}
We need the the expressions of $\theta_j(2^\ell\tau)$ $(j=2,3,4,\ell=1,2,3)$ 
in terms of $\theta_j:=\theta_j(\tau)$ $(j=2,3,4)$, which will be used to prove 
Theorems \ref{Thm2}, \ref{Thm4}, and \ref{Thm3}. 
Recall that the following identities hold for 
any $\tau\in\mathbb{C}$  with $\Im(\tau)>0$:
\begin{eqnarray}
2\theta_2^2(2\tau) &=&\theta_3^2 - \theta_4^2
\,,
\label{Formula1} \\
2\theta_3^2(2\tau) &=&  \theta_3^2 + \theta_4^2
\,,
\label{Formula2}\\
\theta_4^2(2\tau) &=& \theta_3
\theta_4
 \,,
\label{Formula3}
\end{eqnarray} 
and 
\begin{eqnarray}
2\theta_2(4\tau) &=&  \theta_3 - \theta_4\,,
\label{Formula4} \\
2\theta_3(4\tau) &=&  \theta_3+ \theta_4 \,,
\label{Formula5}\\
2\theta_4^4(4\tau) &=&{\big( \theta_3^2+ \theta_4^2\big)}
\theta_3\theta_4\,.
\label{Formula6}
\end{eqnarray}
The most important tool to transfer the algebraic independence of a set of $m$ numbers to another set of $m$ 
numbers, which all satisfy a system of algebraic identities, is given by the following lemma. We call it an
{\em algebraic independence criterion\/} (AIC).
\begin{lemma} \cite[Lemma~3.1]{Elsner}
Let \(x_1,\dots,x_m\in{\C}\) be algebraically independent over \({\Q} \) and let \(y_1,\dots,y_m\in{\C} \) 
satisfy the system of equations 
\[f_j(x_1,\dots,x_m,y_1,\dots,y_m) \,=\, 0 \qquad (1\leq j\leq m) \,,\]
where \(f_j(t_1,\dots,t_m,u_1,\dots,u_m)\in{\Q}[t_1,\dots,t_m,u_1,\dots,u_m] \) \((1\leq j\leq m)\). 
Assume that 
\[\det \left( \frac{\partial f_j}{\partial t_i}(x_1,\dots,x_m,y_1,\dots,y_m) \right) \,\not= \,0 \,.\]
Then the numbers \(y_1,\dots,y_m\) are algebraically independent over \({\Q} \).
\label{Lem1}
\end{lemma}   
We shall apply the AIC to the sets 
$\{ x_1,x_2\}$ with $x_1,x_2\in\Bbb{Z}[\theta_2,\theta_3,\theta_4]$. 
\subsubsection{The case $n=2m$ with odd integer $m$}
In this subsection, 
we put $n=2m$ with an odd integer $m>1$.
In Lemmas~\ref{Lem3.1}, \ref{Lem4}, and \ref{Lem3.2} below, we give sufficient conditions for the numbers in each of the set 
$\{\theta_j(n\tau),\theta_j(\tau)\}$ $(j=2,3,4)$ to be 
algebraically independent over ${\Q}$. 
Replacing $\tau$ by $2\tau$ in (\ref{eq:101}), (\ref{eq:102}), and (\ref{eq:103}), we have 
$$P_m(X_0,Y_0)=0$$ for 
\begin{eqnarray}
X_0=h_2(2\tau)=\displaystyle m^2\frac{\theta_2^4(n\tau)}{\theta_2^4(2\tau)}
&\quad \mbox{and}\quad &
Y_0=
16\frac{\lambda(2\tau)-1}{\lambda(2\tau)}=
16\frac{\theta_2^4(2\tau)-\theta_3^4(2\tau)}{\theta_2^4(2\tau)}, 
\label{eq:160218} \\
X_0=h_3(2\tau)=m^2\frac{\theta_3^4(n\tau)}{\theta_3^4(2\tau)}
&\mbox{and}
&
Y_0=
16\lambda(2\tau)=16\frac{\theta_2^4(2\tau)}{\theta_3^4(2\tau)},\label{eq:160219}\\
X_0=h_4(2\tau)=\, m^2\frac{\theta_4^4(n\tau)}{\theta_4^4(2\tau)} 
&\mbox{and}
&
Y_0=16\frac{\lambda(2\tau)}{\lambda(2\tau)-1}=16\frac{\theta_2^4(2\tau)}{\theta_2^4(2\tau)-\theta_3^4(2\tau)},\label{eq:160220}
\end{eqnarray}
respectively. \\

Let $q$ be an algebraic number with \(q=e^{\pi i \tau} \) and \(\Im (\tau) >0 \). 
The total degree of $P_m(X,Y)$ is denoted by $M$.  
\begin{lemma}
If the polynomial
\[\mbox{\rm Res\,}_X \left( P_m(X,Y),\,X\frac{\partial}{\partial X} P_m\big( X,Y\big) + 
2\big( Y-16 \big)\frac{\partial}{\partial Y} P_m\big( X,Y\big)
\right) \]
does not vanish identically, then the numbers $\theta_2(n\tau)$ and $\theta_2(\tau)$ are algebraically independent over ${\Q}$.
\label{Lem3.1}
\end{lemma}
{\em Proof of Lemma~\ref{Lem3.1}.\/} 
Let 
\[\begin{array}{lcllcl}
x_1 &:=& (\theta_3^4-\theta_4^4)^2\,,\quad & x_2&:=& (\theta_3^2+\theta_4^2)^2 \,,\\
y_1 &:=& 4m^2\theta_2^4(n\tau)\,,\quad & y_2&:=& \theta_2^4 \,.
\end{array} \]     
Then the numbers $x_1$ and $x_2$ are  algebraically independent over ${\Q}$. 
Indeed, the numbers $\theta_3$ and $\theta_4$ are the roots of polynomial
$$
T^8-\frac{1}{2}\left(\frac{x_1}{x_2}+x_2\right)T^4+\frac{1}{16}
\left(\frac{x_1}{x_2}-x_2\right)^2,
$$
so that the field $E:=\mathbb{Q}(\theta_3,\theta_4)$ is an algebraic extension  
of $F:=\mathbb{Q}(x_1,x_2)$, and hence by Lemma~\ref{Lem2}
$$
\begin{array}{ll}
{\rm trans.}\deg F/\mathbb{Q}={\rm trans.}\deg E/F+
{\rm trans.}\deg F/{\mathbb{Q}}=
{\rm trans.}\deg E/{\mathbb{Q}}=2.
\end{array}
$$

By the identities (\ref{Formula1}) and (\ref{Formula2}), the numbers $X_0$ and $Y_0$ 
in (\ref{eq:160218}) are expressed as 
$$
X_0 \,=\, \frac{x_2y_1}{{x_1}} \qquad \mbox{and} \qquad 
Y_0 \,=\, 16\frac{x_1-x_2^2}{{x_1}}
\,.
$$
Define  
\begin{eqnarray}
g_1(t_1,t_2,u_1,u_2)&:=&\frac{t_2u_1}{{t_1}},\nonumber \\\nonumber \\
g_2(t_1,t_2,u_1,u_2)&:=&16\frac{t_1-t_2^2}{{t_1}},\nonumber \\\nonumber \\
f_1(t_1,t_2,u_1,u_2) &:=&t_1^{M} P_m
(g_1,g_2)  \label{H400} \\\nonumber \\
f_2(t_1,t_2,u_1,u_2) &:=& u_2^2 - t_1 \nonumber \,.
\end{eqnarray} 
Since $g_1(x_1,x_2,y_1,y_2)=X_0$ and $g_1(x_1,x_2,y_1,y_2)=Y_0$, 
$$
f_1(x_1,x_2,y_1,y_2)=x_1^{M}P_m(X_0,Y_0)=0
$$
and by the identity (2.1) 
$$
f_2(x_1,x_2,y_1,y_2)=y_2^2-x_1=0.
$$
Using the algebraic independence criterion we have to show the nonvanishing of
\[\Delta \,:=\, \det \left( \begin{array}{cc}
\displaystyle \frac{\partial f_1}{\partial t_1} & \displaystyle \frac{\partial f_1}{\partial t_2} \\ \\
\displaystyle \frac{\partial f_2}{\partial t_1} & \displaystyle \frac{\partial f_2}{\partial t_2}
\end{array} \right) \,=\, \frac{\partial f_1}{\partial t_2} \]
at $(x_1,x_2,y_1,y_2)$; namely by  (\ref{H400})
\begin{eqnarray*}
\frac{\partial f_1}{\partial t_2} \big( x_1,x_2,y_1,y_2 \big)&=& x_1^{M} \frac{\partial P_m}{\partial t_2}(X_0,Y_0)\neq0.
\end{eqnarray*}  
Applying the chain rule, we obtain
\begin{eqnarray*}
\frac{\partial P_m}{\partial t_2}(X_0,Y_0)&=&\frac{\partial P_m}{\partial X}(X_0,Y_0)\cdot\frac{\partial g_1}{\partial t_2}(x_1,x_2,y_1,y_2)+
\frac{\partial P_m}{\partial Y}(X_0,Y_0)\cdot\frac{\partial g_2}{\partial t_2}(x_1,x_2,y_1,y_2)\\\\
&=& \frac{1}{x_2}\Big( \,X_0\frac{\partial P_m}{\partial X}(X_0,Y_0) + 2\big( Y_0-16 \,\big)  
\frac{\partial P_m}{\partial Y}(X_0,Y_0) \,\Big) \,.
\end{eqnarray*}     
Therefore, in order to prove the lemma by the algebraic independence criterion, 
it suffices to show that
\begin{equation}
X_0\frac{\partial P_m}{\partial X}(X_0,Y_0) + 2\big( Y_0 -16 \,\big) \frac{\partial P_m}{\partial Y}(X_0,Y_0) \,\not= \, 0\,.
\label{H500}            
\end{equation}        

By the hypothesis of Lemma~\ref{Lem3.1} the polynomial
\[R(Y):=\mbox{Res\,}_X \left( P_m(X,Y),\,X\frac{\partial}{\partial X} P_m\big( X,Y\big) + 
2\big( Y-16 \big)\frac{\partial}{\partial Y} P_m\big( X,Y\big)
\right)\in{\mathbb{Z}}[Y] \]
does not vanish identically. 
For fixed $Y=Y_0:=16(x_1-x_2^2)/x_1$ we have $R(Y_0)\in {\Q}(x_1,x_2)$, so that the algebraic independence of $x_1,x_2$ proves
$R(Y_0)\not= 0$. 
In particular, $P_m(X,Y_0)$ and 
$$
X\frac{\partial}{\partial X} P_m\big( X,Y_0\big) + 
2\big( Y_0-16 \big)\frac{\partial}{\partial Y} P_m\big( X,Y_0\big)
$$
(which both are polynomials in $X$) have no common root. 
Since $P_m(X,Y_0)$ vanishes for $X=X_0:=x_2y_1/x_1$, we obtain (\ref{H500}).         
The proof of Lemma~\ref{Lem3.1} is completed.   
\hfill $\Box$                 
\par\medskip      

\begin{lemma}
If the polynomial
\[\mbox{\rm Res\,}_X \left( P_m(X^2,Y^2),\,X^2\frac{\partial}{\partial X} P_m\big( X^2,Y^2\big) + \big( Y^2+4Y \big)\frac{\partial}{\partial Y} P_m\big( X^2,Y^2\big)
\right) \]
does not vanish identically, then the numbers $\theta_3(n\tau)$ and $\theta_3(\tau)$ are algebraically independent over ${\Q}$.
\label{Lem4}
\end{lemma}
{\em Proof of Lemma~\ref{Lem4}.\/} 
Let
\[\begin{array}{lcllcl}
x_1 &:=& 2\theta_3^2\,,\quad & x_2 &:=& \theta_3^2+\theta_4^2 \,,\\
y_1 &:=& 2m\theta_3^2(n\tau)\,,\quad & y_2 &:=&\, \theta_3^2 \,.
\end{array} \]      
Then the numbers $x_1,x_2$ are algebraically independent over ${\Q}$ and 
we see by (\ref{Formula1}) and (\ref{Formula2}) that the numbers 
$X_0$ and $Y_0$ in (\ref{eq:160219}) are given by 
\begin{equation}
X_0 \,=\, \frac{y_1^2}{x_2^2} \qquad \mbox{and} \qquad 
Y_0=(\sqrt{Y_0})^2 \,:=\, \left(\frac{4{(x_1 - x_2)}}{x_2}\right)^2 \,.
\label{H300}
\end{equation}                  
Define 
\begin{eqnarray}
f_1(t_1,t_2,u_1,u_2) &:=&{t_2}^{2M} P_m
\left(\frac{u_1^2}{t_2^2},\frac{16{(t_1 - t_2)}^2}{t_2^2}\right) \nonumber \\\nonumber \\
f_2(t_1,t_2,u_1,u_2) &:=& 2u_2 - t_1. \nonumber \\\nonumber \,
\end{eqnarray} 
Similarly to the proof of Lemma \ref{Lem3.1}, 
applying the algebraic independence criterion, we have to show that 
\begin{equation}
\frac{\partial P_m}{\partial t_2}(X_0,Y_0)
=-\frac{2}{x_2}\Big( \,X_0\frac{\partial P_m}{\partial X}(X_0,Y_0) + \big( Y_0 + 4\sqrt{Y_0} \,\big)  
\frac{\partial P_m}{\partial Y}(X_0,Y_0) \,\Big)\neq0.
\label{eq:160221}
\end{equation}     
By the hypothesis of Lemma~\ref{Lem4} the polynomial
\[R(Y) \,:=\, \mbox{Res\,}_X \left( P_m(X^2,Y^2),\,X^2\frac{\partial P_m}{\partial X}\big( X^2,Y^2\big) + \big( Y^2+4Y \big)\frac{\partial P_m}{\partial Y}
\big( X^2,Y^2\big) \right) \,\in \, {\Z}[Y] \]
does not vanish identically. 
Since the numbers $x_1$ and $x_2$ are algebraically independent over ${\Q}$, 
we have $R(Y_1)\not= 0$ for $Y_1:=4(x_1-x_2)/x_2$, and hence 
the polynomials $P_m(X^2,Y_1^2)$ and 
\[X^2\frac{\partial P_m}{\partial X}\big( X^2,Y_1^2\big) + \big( Y_1^2+4Y_1 \big)\frac{\partial P_m}{\partial Y} \big( X^2,Y_1^2\big) \]
have no common root. 
Noting that $P_m(X_1^2,Y_1^2)=0$ holds for $X_1:=y_1/x_2$, we obtain               
\[X_1^2\frac{\partial P_m}{\partial X}\big( X_1^2,Y_1^2\big) + \big( Y_1^2+4Y_1 \big)\frac{\partial P_m}{\partial Y} \big( X_1^2,Y_1^2\big) \,\not=\, 0\,.\]                
Finally, using $X_0=X_1^2$ and $Y_0=Y_1^2$ by (\ref{H300}), we obtain (\ref{eq:160221}), which completes the proof of Lemma~\ref{Lem4}. \hfill $\Box$                 
\par\medskip      
\begin{lemma}
If the polynomial
$$
\mbox{\rm Res\,}_X 
\left( P_m(X,Y),\,
X^2\left(
\frac{\partial P_m}{\partial X}(X,Y)\right)^2-
Y\big( Y -16 \,\big) 
\left(\frac{\partial P_m}{\partial Y}(X,Y)\right)^2
\right)
$$
does not vanish identically, then the numbers $\theta_4(n\tau)$ and $\theta_4(\tau)$ are algebraically independent over ${\Q}$.
\label{Lem3.2}
\end{lemma}
{\em Proof of Lemma~\ref{Lem3.2}.\/} 
By (\ref{Formula1}), (\ref{Formula2}), and (\ref{Formula3}), the numbers 
$X_0$ and $Y_0$ in (\ref{eq:160220}) are expressed as 
$$
X_0 \,=\, \frac{y_1}{{x_2}y_2} \qquad \mbox{and} \qquad 
Y_0 \,=\, -4\frac{(x_2-y_2)^2}{{x_2}y_2}
$$   
with 
$$
Y_0(Y_0-16)=(\sqrt{Y_0(Y_0-16)})^2:=
\left(
\frac{4(x_2^2-y_2^2)}{{x_2}y_2}
\right)^2
\,,
$$
where 
\[\begin{array}{lcllcl}
x_1 &:=&\theta_4^2 \,,\quad & x_2&:=&\theta_3^2 \,,\\
y_1 &:=& m^2\theta_4^4(n\tau)\,,\quad & y_2&:=& \theta_4^2\,.
\end{array} \]  
Define 
\begin{eqnarray}
f_1(t_1,t_2,u_1,u_2) &:=&(t_2u_2)^{M} P_m
\left(\frac{u_1}{{t_2}u_2},-4\frac{(t_2-u_2)^2}{{t_2}u_2}\right)  \nonumber \\\nonumber \\
f_2(t_1,t_2,u_1,u_2) &:=& u_2 - t_1 \nonumber \,.
\end{eqnarray} 
Then, similarly to the proofs of privious lemmas, we have only to prove
\begin{equation}
\frac{\partial P_m}{\partial t_2}(X_0,Y_0)
=-\frac{1}{x_2}\Big( \,X_0\frac{\partial P_m}{\partial X}(X_0,Y_0) -
\sqrt{Y_0\big( Y_0-16 \,\big)}  
\frac{\partial P_m}{\partial Y}(X_0,Y_0) \,\Big)\neq0 \,,
\end{equation}     
which follows immediately from the hypothesis of Lemma~\ref{Lem3.2}. 
\hfill $\Box$                 
\subsubsection{The case $n=4m$ with odd integer $m$}
\begin{lemma}
Let $n=4m$, where $m>1$ is an odd integer. Let $q$ be an algebraic number with $q=e^{\pi i\tau}$ and $\Im (\tau)>0$. If the polynomial
\[\mbox{Res\,}_X \left( P_m(X^4,Y^4),\,X^4\frac{\partial}{\partial X} P_m\big( X^4,Y^4\big) + \big( Y^4+2Y^3 \big)\frac{\partial}{\partial Y} P_m\big( X^4,Y^4\big)
\right) \]
does not vanish identically, then the numbers $\theta_3(n\tau)$ and $\theta_3(\tau)$ are algebraically independent over ${\Q}$.
\label{Lem5}
\end{lemma}
{\em Proof of Lemma~\ref{Lem5}.\/} 
Let 
\[\begin{array}{lcllcl}
x_1 &:=&2\theta_3 \,,\quad & x_2&:=&\theta_3+\theta_4 \,,\\
y_1 &:=& 16m^2\theta_3^4(n\tau)\,,\quad & y_2&:=& \theta_3 \,.
\end{array} \]  
Then, by the identities (\ref{Formula4}) and (\ref{Formula5}), the polynomial 
$P_m(X,Y)$ vanishes at
$$
X_0 \,=\, \frac{y_1}{x_2^4} \qquad \mbox{and} \qquad 
Y_0=(\sqrt[4]{Y_0})^4 \,:=\, \left(\frac{2{(x_1 - x_2)}}{x_2}\right)^4\quad \mbox{with}\quad 
\sqrt[4]{Y_0^3}:=(\sqrt[4]{Y_0})^3.
$$                  
Again we have $2y_2-x_1=0$, and $x_1,x_2$ are algebraically independent over ${\Q}$ for any algebraic number $q=e^{\pi i\tau}$ with $\Im (\tau)>0$. We introduce 
the polynomials
\begin{eqnarray*}
f_1(t_1,t_2,u_1,u_2) &:=&t_2^{4M} P_m\Big( \,\frac{u_1}{t_2^4},
\frac{16{(t_1-t_2)}^4}{t_2^4}\,\Big) \,, 
\label{H100} \\
f_2(t_1,t_2,u_1,u_2) &:=& 2u_2 - t_1 \nonumber \,.
\end{eqnarray*} 
Using the algebraic independence criterion, we have to show that 
\begin{equation}
\frac{\partial f_1}{\partial t_2} \big( x_1,x_2,y_1,y_2 \big)=-4x_2^{4M-1} \Big( \,X_0\frac{\partial P_m}{\partial X}(X_0,Y_0) + \big( Y_0 + 2\sqrt[4]{Y_0^3} \,\big)  
\frac{\partial P_m}{\partial Y}(X_0,Y_0) \,\Big)\neq0 \,.
\end{equation}      
which follows from the hypothesis of Lemma~\ref{Lem5}. \hfill $\Box$   
\par\medskip
\section{Proof of Theorems} \label{Sec4}  

{\it Proof of Theorem \ref{Thm1}}. 
We consider the case of  $\theta_3(\tau)$. Let  
$$
\begin{array}{ll}
\alpha_1:=16\theta_2^4,&\quad \alpha_2:=\theta_3^4,\\\\
\beta_1:=n^2\theta_3^4(n\tau),&\quad \beta_2:=\theta_3^4,
\end{array}
$$ 
where, by Lemmas~\ref{Lem2} and \ref{Lem7}, 
the numbers $\alpha_1$ and $\alpha_2$ are algebraically independent over 
$\Bbb{Q}$ and the number $\beta_1/\beta_2=h_3(\tau)$ is transcendental.
Define $P(X,Y):=P_n(Y,X)$ and $Q(X,Y):=Y$. By 
(\ref{eq:102})
$$
P(\alpha_1/\alpha_2,\beta_1/\beta_2)=P_n(h_3(\tau),16\lambda(\tau))=0
$$
and 
$$
Q(\beta_1,\beta_2)=\beta_2=\alpha_2.
$$
Hence, applying Lemma \ref{Lem71}, we obtain the algebraic independence over $\mathbb{Q}$ of the numbers $\beta_1$ and $\beta_2$. This implies that the numbers 
$\theta_3(n\tau)$ and $\theta_3(\tau)$ are algebraically independent over $\mathbb{Q}$.\\
The same holds for the sets $\{\theta_2(n\tau),\theta_2(\tau)\}$ and $\{\theta_4(n\tau),\theta_4(\tau)\}$. In these cases, we use the identities (\ref{eq:101}), (\ref{eq:103}) and Lemma \ref{Lem71} with 
$$
\begin{array}{ll}
\alpha_1:=16(\theta_2^4-\theta_3^4),&\quad \alpha_2:=\theta_2^4,\\\\
\beta_1:=n^2\theta_2^4(n\tau),&\quad \beta_2:=\theta_2^4,\\\\
P(X,Y):=P_n(Y,X),&\quad Q(X,Y):=Y,
\end{array}
$$
and 
$$
\begin{array}{ll}
\alpha_1:=16\theta_2^4,&\quad \alpha_2:=\theta_2^4-\theta_3^4,\\\\
\beta_1:=n^2\theta_4^4(n\tau),&\quad \beta_2:=\theta_4^4,\\\\
P(X,Y):=P_n(Y,X),&\quad Q(X,Y)=-Y,
\end{array}
$$
respectively. In the latter case, we note that the equality 
$Q(\beta_1,\beta_2)=\alpha_2$ holds from 
the identity (\ref{eq:405}). 
Thus, the proof of Theorem \ref{Thm1} is completed.
\hfill \qed\\

{\it Proof of Theorem~\ref{Thm2}.} 
We first consider the case $n=2$. 
Let $F:=\mathbb{Q}(x_1,x_2)$, where $x_1:=2\theta_2^2(2\tau)$ and 
$x_2:=\theta_2^4$. 
Then by the identity (\ref{Formula1}) together with the relation 
$\theta_2^4=\theta_3^4-\theta_4^4$, we have $F\subset E:=\mathbb{Q}(\theta_3,\theta_4)$ and 
$$
2x_1\theta_3^2-x_1^2-x_2=2x_1\theta_4^2+x_1^2-x_2=0.
$$
This implies that the field extension $E/F$ is algebraic, so that 
$$
{\rm trans}\deg F/\mathbb{Q}=
{\rm trans}\deg E/F+ {\rm trans}\deg F/\mathbb{Q}
={\rm trans}\deg E/\mathbb{Q}=2,
$$
which implies that the 
numbers $\theta_2(2\tau)$ and $\theta_2(\tau)$ are 
are algebraically independent 
over $\mathbb{Q}$. 
For $n=4$, putting $F:=\mathbb{Q}(2\theta_2(4\tau),\theta_2^4)$
and using (\ref{Formula4}), we can proceed the same argument as stated above. 
For $x_1:=2\theta_2(4\tau)$ and $x_2:=\theta_2^4$ we use the identities
\[\theta_3^4-{\big( x_1-\theta_3 \big)}^4 - x_2 \,=\, \theta_4^4-{\big( x_1+\theta_4 \big)}^4 + x_2 \,=\, 0\,. \]
In the case of $n=6$, we use Lemma~\ref{Lem3.1}. 
From \cite{Nes3} we know that
\[P_3(X,Y) \,=\, 9-\big( Y^2-16Y+28\big) X+30X^2-12X^3+X^4 \,.\]
Hence, we have
\begin{eqnarray*}
&& \mbox{\rm Res\,}_X \left( P_3(X,Y),\,X\frac{\partial}{\partial X} P_3\big( X,Y\big) + 
2\big( Y-16 \big)\frac{\partial}{\partial Y} P_3\big( X,Y\big) \right) \\
&=& 9Y\big( 5Y-512\big) {\big( Y-16\big)}^2{\big( Y-8\big)}^4 
\,\not\equiv \,0 \,.
\end{eqnarray*}
Lemma~\ref{Lem3.1} gives the desired result for $n=6$. \hfill \qed\\

{\it Proof of Theorem~\ref{Thm4}.} 
Similarly to the proof of Theorem~\ref{Thm2}, 
we can deduce the conclusion for the cases $n=2,4,8$ by using the identities 
(\ref{Formula3}), (\ref{Formula6}), and 
$$
32\theta_4^8(8\tau)=
\left(\theta_3+\theta_4\right)^4
\left(\theta_3^2+\theta_4^2\right)\theta_3\theta_4,
$$
which is yielded from (\ref{Formula2}), (\ref{Formula3}), and (\ref{Formula6}). 
In the cases of $n=6,10$, we use Lemma~\ref{Lem3.2}. For $n=6$ we compute the resultant from the lemma explicitly.
\begin{eqnarray*}
&& \mbox{\rm Res\,}_X \left( P_3(X,Y),\,X^2\left( \frac{\partial P_3}{\partial X}(X,Y)\right)^2-Y\big( Y -16 \,\big) 
\left(\frac{\partial P_3}{\partial Y}(X,Y)\right)^2 \right) \\
&=& -81Y^3\big( 375Y^2-6000Y+262144\big) {\big( Y-16\big)}^3{\big( Y-8\big)}^8 \,\not\equiv \,0 \,.
\end{eqnarray*}
Next, let $n=10$. In \cite{Nes3} the polynomial $P_5(X,Y)$ is given as well.
\begin{eqnarray*}
P_5(X,Y) &\,=\,& 25 - (126 - 832Y + 308Y^2 - 32Y^3 + Y^4)X + (255 + 1920Y - 120Y^2)X^2 \\
&& +\,(-260 + 320Y - 20Y^2)X^3 + 135X^4 - 30X^5 + X^6 \,.
\end{eqnarray*}
Hence, by setting
\[T_{10}(Y) \,:=\, \mbox{\rm Res\,}_X \left( P_5(X,Y),\,X^2\left( \frac{\partial P_5}{\partial X}(X,Y)\right)^2-Y\big( Y -16 \,\big) 
\left(\frac{\partial P_5}{\partial Y}(X,Y)\right)^2 \right) \,,\]
we obtain $T_{10}(1) \equiv 1 \pmod 2$, such that $T_{10}(Y)$ does not vanish identically. This completes the proof of Theorem~\ref{Thm4}. \hfill \qed\\


{\it Proof of Theorem~\ref{Thm3}.} 
Taking the results from Theorem~A and Theorem~\ref{Thm1} into account, for Theorem~\ref{Thm3} it suffices to consider $n\in \{ 14,18,20,22\}$. 
Here, we compute the resultants from Lemma~\ref{Lem4} (for $n\in \{ 14,18,22\}$) and from Lemma~\ref{Lem5} (for $n=20$) explicitly by using a computer algebra system.
In order to show that the resultants do not vanish we again consider the values at $Y=1$. For $n=14$ we use
\begin{eqnarray*}
P_7(X,Y) &\,=\,& 49 - (344 - 17568Y + 20554Y^2 - 6528Y^3 + 844Y^4 - 48Y^5 + Y^6)X \\
&& +\,(1036 + 156800Y + 88760Y^2 - 12320Y^3 + 385Y^4)X^2 \\
&& -\,(1736 - 185024Y + 18732Y^2 - 896Y^3 + 28Y^4)X^3 \\
&& +\,(1750 + 31360Y - 1960Y^2)X^4 - (1064 - 2464Y + 154Y^2)X^5 \\
&& +\,364X^6 - 56X^7 + X^8
\end{eqnarray*}
(cf. \cite{Elsner3}) to obtain the resultant from Lemma~\ref{Lem4},
\[T_{14}(Y) \,:=\, \mbox{Res\,}_X \left( P_7(X^2,Y^2),\,X^2\frac{\partial}{\partial X} P_7\big( X^2,Y^2\big) + \big( Y^2+4Y \big)\frac{\partial}{\partial Y} 
P_7\big( X^2,Y^2\big) \right) \,.\]
It follows that $T_{14}(1) \equiv 1 \pmod 2$. This shows that $T_{14}(1)\not= 0$.
For $n=18$ we apply Lemma~\ref{Lem4} with
\begin{eqnarray*}
P_9(X,Y) &\,=\,& 6561 - (60588-18652032Y+56033208Y^2-40036032Y^3+11743542Y^4 \\
&& -\,1715904Y^5+132516Y^6-5184Y^7+81Y^8)X \\
&& +\,(250146+427613184Y+2083563072Y^2+86274432Y^3-57982860Y^4 \\
&& +\,4249728Y^5-99288Y^6+576Y^7-9Y^8)X^2 \\
&& -\,(607420-1418904064Y+2511615520Y^2-353755456Y^3+19071754Y^4 \\
&& -\,612736Y^5+13960Y^6-64Y^7+Y^8)X^3 \\
&& +\,(959535+856286208Y+8468928Y^2-2145024Y^3-808488Y^4 \\
&& +\,65664Y^5-1368Y^6)X^4 \\
&& -\,(1028952+22899456Y+1430352Y^2-505152Y^3+38826Y^4 \\
&& -\,1728Y^5+36Y^6)X^5 \\
&& +\,(757596-13138944Y+4160448Y^2-417408Y^3+13044Y^4)X^6 \\
&& -\,(378072+1138176Y+16416Y^2-10944Y^3+342Y^4)X^7 \\
&& +\,(122895+64512Y-4032Y^2)X^8 - (24060-11136Y+696Y^2)X^9 \\
&& +\,2466X^{10} - 108X^{11} + X^{12} \,.
\end{eqnarray*}
We have
\[T_{18}(Y) \,:=\, \mbox{Res\,}_X \left( P_9(X^2,Y^2),\,X^2\frac{\partial}{\partial X} P_9\big( X^2,Y^2\big) + \big( Y^2+4Y \big)\frac{\partial}{\partial Y} 
P_9\big( X^2,Y^2\big) \right) \,,\]
and thus $T_{18}(1) \equiv 1 \pmod 2$. Hence, $T_{18}(1)\not= 0$.
For $n=20$ we need the polynomial $P_5(X,Y)$, which was already used in the proof of Theorem~\ref{Thm4}. 
The resultant from Lemma~\ref{Lem5},
\[T_{20}(Y) \,:=\, \mbox{Res\,}_X \left( P_5(X^4,Y^4),\,X^4\frac{\partial}{\partial X} P_5\big( X^4,Y^4\big) + \big( Y^4+2Y^3 \big)\frac{\partial}{\partial Y} 
P_5\big( X^4,Y^4\big) \right) \,,\]
satisfies $T_{20}(1) \equiv 1 \pmod 2$.
Finally, for $n=22$ we again apply Lemma~\ref{Lem4} with
\begin{eqnarray*}
P_{11}(X,Y) &\,=\,& 121 - (1332-2214576Y+15234219Y^2-21424896Y^3+11848792Y^4 \\
&& -\,3309152Y^5+522914Y^6-48896Y^7+2684Y^8-80Y^9+Y^{10})X \\
&& +\,(6666+111458688Y+2532888424Y^2+2367855776Y^3-327773413Y^4 \\
&& -\,9982720Y^5+3230480Y^6-161920Y^7+2530Y^8)X^2 \\
&& -\,(20020-864654912Y+12880909668Y^2-5289254784Y^3+744094076Y^4 \\
&& -\,43914992Y^5+967461Y^6-2816Y^7+44Y^8 )X^3 \\
\end{eqnarray*}
\begin{eqnarray*}
&& +\,(40095+1748954240Y-175142088Y^2+372281536Y^3-68516998Y^4 \\
&& +\,4266240Y^5-88880Y^6)X^4 \\
&& -\,(56232-1061669664Y+132688050Y^2-10724736Y^3+715308Y^4 \\
&& -\,28512Y^5+594Y^6)X^5 \\
&& +\,(56364+211953280Y-7454568Y^2-724064Y^3+22627Y^4)X^6 \\
&& -\,(40392-24140864Y+2162116Y^2-81664Y^3+2552Y^4)X^7 \\
&& +\,(20295+1448832Y-90552Y^2)X^8 - (6820-36784Y+2299Y^2)X^9 \\
&& +\,1386X^{10} - 132X^{11} + X^{12} \,.
\end{eqnarray*}
Here, we obtain
\[T_{22}(Y) \,:=\, \mbox{Res\,}_X \left( P_{11}(X^2,Y^2),\,X^2\frac{\partial}{\partial X} P_{11}\big( X^2,Y^2\big) + \big( Y^2+4Y \big)\frac{\partial}{\partial Y} 
P_{11}\big( X^2,Y^2\big) \right) \,,\]
where $T_{22}(1) \equiv 3 \pmod {13}$, such that $T_{22}(1) \not= 0$. 
The proof of Theorem~\ref{Thm2} is complete. \hfill \qed\\

{\bf Acknowledgements.}
The second author was supported by Japan Society for the Promotion of Science, Grant-in-Aid for Young Scientists (B), 15K17504.


\end{document}